\documentclass{article}
\usepackage{amssymb}
\usepackage{amsfonts}
\usepackage{amsmath}

\input{tcilatex}

\begin{document}

\begin{center}
\textbf{Latter research on Euler-Mascheroni constant}%
\begin{equation*}
\end{equation*}

Valentin Gabriel Cristea$^{1}$ and Cristinel Mortici$^{2}$%
\begin{equation*}
\end{equation*}

$^{1}$Ph. D. Student, University Politehnica of Bucharest, Splaiul Independen%
\c{t}ei 313, Bucharest, Romania, valentingabrielc@yahoo.com\bigskip

$^{2}$Prof. dr. habil., Valahia University of T\^{a}rgovi\c{s}te, Bd. Unirii
18, 130082 T\^{a}rgovi\c{s}te, Romania, cristinelmortici@yahoo.com%
\begin{equation*}
\end{equation*}
\end{center}

\begin{quotation}
\textbf{Abstract: }In this work, we present a review and an example on some
latter results on the problem of approximating the Euler-Mascheroni
constant. We use the method firstly introduced in [C. Mortici, Product
Approximations via Asymptotic Integration Amer. Math. Monthly 117 (5) (2010)
434-441].
\end{quotation}

\bigskip

\begin{equation*}
\end{equation*}

\textbf{1. Introduction and motivation}\bigskip

The famous Euler-Mascheroni constant $\gamma =0,57721566490153286...$ was
firstly studied by the Swiss mathematician Leonhard Euler (1707-1783) and
the Italian mathematician Lorenzo Mascheroni (1750-1800) as the limit of the
sequence%
\begin{equation*}
\gamma _{n}=1+\frac{1}{2}+\frac{1}{3}+\cdots +\frac{1}{n}-\ln n.
\end{equation*}%
The question on $\gamma $ whether it is a rational number or not, has not an
answer yet. The reason seems to be the inexistence of very fast convergences
to $\gamma $, having a simple form. As a consequence, in the recent past,
many authors introduced new fast convergences to $\gamma .$

We start a history of such approximations with the results of Tims and
Tyrrell [20]%
\begin{equation*}
\frac{1}{2(n+1)}<\gamma _{n}-\gamma <\frac{1}{2(n-1)},
\end{equation*}%
Young [23]%
\begin{equation*}
\frac{1}{2(n+1)}<\gamma _{n}-\gamma <\frac{1}{2n},
\end{equation*}%
Anderson, Barnard, Richards, Vamanamurthy and Vuorinen [3]%
\begin{equation*}
\frac{1-\gamma }{n}<\gamma _{n}-\gamma <\frac{1}{2n},
\end{equation*}%
Mortici and Vernescu [13]%
\begin{equation*}
\frac{1}{2n+1}<\gamma _{n}-\gamma <\frac{1}{2n},
\end{equation*}%
or Toth [21]%
\begin{equation*}
\frac{1}{2n+\frac{2}{5}}<\gamma _{n}-\gamma <\frac{1}{2n+\frac{1}{3}}.
\end{equation*}%
Independently, Alzer [2] and Chen-Qi [4] proved%
\begin{equation*}
\frac{1}{2n+\frac{2\gamma -1}{1-\gamma }}<\gamma _{n}-\gamma <\frac{1}{2n+%
\frac{1}{3}}.
\end{equation*}%
Next, Qiu and Vuorinen [19, Cor. 2.13] showed%
\begin{equation*}
\frac{1}{2n}-\frac{\alpha }{n^{2}}<\gamma _{n}-\gamma <\frac{1}{2n}-\frac{%
\beta }{n^{2}},
\end{equation*}%
where $\alpha =1/2$ and $\beta =\gamma -1/2,$\ with its consequence%
\begin{equation*}
\frac{1}{2n}-\frac{1}{8n^{2}}<\gamma _{n}-\gamma <\frac{1}{2n}.
\end{equation*}%
This is also called Franel's inequality [18, Ex. 18]. Karatsuba [9] proved%
\begin{equation*}
\frac{1}{2n}-\frac{1}{12n^{2}}+\frac{1}{120n^{4}}-\frac{1}{126n^{6}}<\gamma
_{n}-\gamma <\frac{1}{2n}-\frac{1}{12n^{2}}+\frac{1}{120n^{4}},
\end{equation*}%
while Mortici [14] gave%
\begin{equation*}
\frac{1}{2n+\frac{1}{3}+\frac{1}{18n}}<\gamma _{n}-\gamma <\frac{1}{2n+\frac{%
1}{3}+\frac{1}{32n}}.
\end{equation*}%
DeTemple [7, 8] introduced the sequence%
\begin{equation*}
R_{n}=1+\frac{1}{2}+\frac{1}{3}+\cdots +\frac{1}{n}-\ln \left( n+\frac{1}{2}%
\right) 
\end{equation*}%
that converges to $\gamma $ like $n^{-2}$ , since%
\begin{equation*}
\frac{1}{24(n+1)^{2}}<R_{n}-\gamma <\frac{1}{24n^{2}}.
\end{equation*}%
Chen [5] found%
\begin{equation*}
\frac{1}{24(n+a)^{2}}<R_{n}-\gamma <\frac{1}{24(n+b)^{2}},
\end{equation*}%
with $a=$ $\frac{1}{\sqrt{24\left( -\gamma +1-\ln \frac{3}{2}\right) }}-1$
and $b=\frac{1}{2}$.

Chen and Mortici [6] improved these bounds to%
\begin{eqnarray*}
&&\frac{1}{24(n+\frac{1}{2})^{2}}-\frac{\frac{7}{960}}{\left( n+\frac{1}{2}%
\right) ^{4}}+\frac{\frac{31}{8064}}{\left( n+\frac{1}{2}\right) ^{6}}-\frac{%
\frac{127}{30720}}{\left( n+\frac{1}{2}\right) ^{8}} \\
&<&R_{n}-\gamma \\
&<&\frac{1}{24(n+\frac{1}{2})^{2}}-\frac{\frac{7}{960}}{\left( n+\frac{1}{2}%
\right) ^{4}}+\frac{\frac{31}{8064}}{\left( n+\frac{1}{2}\right) ^{6}}
\end{eqnarray*}%
For further reading, please see [13], [14], [15], [18], [21].

Mortici [15], introduced the sequence%
\begin{equation*}
\mu _{n}\left( a,b\right) =1+\frac{1}{2}+\frac{1}{3}+\cdots +\frac{1}{n-1}+%
\frac{1}{an}-\ln \left( n+b\right) ,
\end{equation*}%
depending on real parameters $a,$ $b.$

This family extends the sequence $V_{n}=\mu _{n}\left( 2,0\right) $
introduced by Vernescu [22] and DeTemple sequence $R_{n}=\mu _{n}\left( 1,%
\frac{1}{2}\right) .$ See [7]. Both sequences $V_{n}$ and $R_{n}$ converges
as $n^{-2}.$ For proofs and details, see [7], [8], [22].

The results in [15] show that for $a=6-2\sqrt{6}$, $b=1/\sqrt{6}$%
\begin{equation*}
u_{n}=1+\frac{1}{2}+\frac{1}{3}+\cdots +\frac{1}{n-1}+\frac{1}{(6-2\sqrt{6})n%
}-\ln \left( n+\frac{1}{\sqrt{6}}\right)
\end{equation*}%
and $a=6+2\sqrt{6}$, $b=-1/\sqrt{6}$%
\begin{equation*}
u_{n}=1+\frac{1}{2}+\frac{1}{3}+\cdots +\frac{1}{n-1}+\frac{1}{(6+2\sqrt{6})n%
}-\ln \left( n-\frac{1}{\sqrt{6}}\right)
\end{equation*}%
converge to $\gamma $ with the speed of convergence at $n^{-3}.$%
\begin{equation*}
\end{equation*}

\textbf{2. A convergence towards }$\gamma $\bigskip 

Using the idea from [15], we introduce the family of sequences $%
v_{n}=v_{n}\left( a,b\right) $%
\begin{equation*}
v_{n}\left( a,b\right) =1+\frac{1}{2}+\cdots +\frac{1}{n-2}+\frac{an+b}{%
n(n-1)}-\ln n,
\end{equation*}%
depending on real parameters $a$ and $b.$ In order to avoid some
inconvenience, we assume $v_{0},$ $v_{1},$ $v_{2}$ given.

This is an extension of the classical convergence $\left( \gamma _{n}\right)
_{n\geq 1},$ since%
\begin{equation*}
\gamma _{n}=v_{n}\left( 2,-1\right) .
\end{equation*}%
Known fact, the new introduced sequence converges to $\gamma $ as $n^{-1}$
in case $a=2,$ $b=-1.$

The problem we rise here is what are the best parameters $a$ and $b$ which
provide the fastest sequence $v_{n}\left( a,b\right) .$ The answer is
formulated as the following\bigskip

\textbf{Theorem 2.1. }(i)\emph{\ If }$a\neq \frac{3}{2},$\emph{\ then the
sequence }$\left( v_{n}\left( a,b\right) \right) _{n\geq 1}$\emph{\ has the
rate of convergence }$n^{-1}$\emph{.}

(ii)\emph{\ If }$a=\frac{3}{2}$\emph{\ and }$b\neq -\frac{5}{12}$\emph{\
then the sequence }$\left( v_{n}\left( \frac{3}{2},b\right) \right) _{n\geq
1}$\emph{\ has the rate of convergence }$n^{-2}$\emph{.}

(iii)\emph{\ If }$a=\frac{3}{2}$\emph{\ and }$b=-\frac{5}{12}$\emph{\ then
the sequence }$\left( v_{n}\left( \frac{3}{2},-\frac{5}{12}\right) \right)
_{n\geq 1}$\emph{\ has the rate of convergence }$n^{-3}$\emph{.}\bigskip

We use the following\bigskip

\textbf{Lemma 2.1. }\emph{If the sequence }$\left( x_{n}\right) _{n\geq 1}$%
\emph{\ is convergent to }$x$ \emph{and there exists the limit}%
\begin{equation*}
\underset{n\rightarrow \infty }{\lim }n^{k}\left( x_{n}-x_{n+1}\right) =l\in 
\mathbb{R}%
\end{equation*}%
\emph{with }$k>1$\emph{, then there exists the limit}%
\begin{equation*}
\underset{n\rightarrow \infty }{\lim }n^{k-1}\left( x_{n}-x\right) =\frac{l}{%
k-1}.
\end{equation*}

This is a form of Cesaro-Stolz lemma, which is useful in constructing of
asymptotic expansions, or evaluating the speed of convergence. For proof and
other details, see \emph{e.g. }[12].\bigskip

\emph{Proof of Theorem 2.1.}\textbf{\ }We have%
\begin{equation*}
v_{n}-v_{n+1}=\frac{an+b}{n(n-1)}-\frac{1}{n-1}-\frac{a\left( n+1\right) +b}{%
n(n+1)}-\ln \frac{n}{n+1}.
\end{equation*}%
By using a computer software such as Maple, we get%
\begin{equation*}
v_{n}-v_{n+1}=\left( a-\frac{3}{2}\right) \frac{1}{n^{2}}+\left( a+2b-\frac{2%
}{3}\right) \frac{1}{n^{3}}
\end{equation*}%
\begin{equation}
+\left( a-\frac{5}{4}\right) \frac{1}{n^{4}}+\left( a+2b-\frac{4}{5}\right) 
\frac{1}{n^{5}}+O\left( \frac{1}{n^{6}}\right) .  \tag{1}
\end{equation}

(i) If $a\neq \frac{3}{2},$ then%
\begin{equation*}
\lim_{n\rightarrow \infty }n^{2}\left( v_{n}-v_{n+1}\right) =a-\frac{3}{2}%
\neq 0,
\end{equation*}%
while Lemma 2.1 says%
\begin{equation*}
\lim_{n\rightarrow \infty }n\left( v_{n}-\gamma \right) =a-\frac{3}{2}\neq 0.
\end{equation*}%
As a consequence, $\left( v_{n}\right) _{n\geq 1}$ converges as $n^{-1}.$

(ii) and (iii). If $a=\frac{3}{2},$ then (1) reads as%
\begin{equation*}
v_{n}-v_{n+1}=\left( 2b+\frac{5}{6}\right) \frac{1}{n^{3}}+\frac{1}{4n^{4}}%
+\left( 2b+\frac{7}{10}\right) \frac{1}{n^{5}}+O\left( \frac{1}{n^{6}}%
\right) .
\end{equation*}%
If $b\neq -\frac{5}{12},$ then%
\begin{equation*}
\lim_{n\rightarrow \infty }n^{3}\left( v_{n}-v_{n+1}\right) =2b+\frac{5}{6}%
\neq 0,
\end{equation*}%
and by Lemma 2.1,%
\begin{equation*}
\lim_{n\rightarrow \infty }n^{2}\left( v_{n}-\gamma \right) =b+\frac{5}{12}%
\neq 0.
\end{equation*}%
As a consequence, $\left( v_{n}\left( \frac{3}{2},b\right) \right) _{n\geq
1} $, with $b\neq -\frac{5}{12},$ converges as $n^{-2}.$

Finally, with $a=\frac{3}{2}$ and $b=-\frac{5}{12},$ we get from (1)%
\begin{equation*}
v_{n}-v_{n+1}=\frac{1}{4n^{4}}-\frac{2}{15n^{5}}+O\left( \frac{1}{n^{6}}%
\right) ,
\end{equation*}%
then use Lemma 2.1 to obtain%
\begin{equation}
\lim_{n\rightarrow \infty }n^{3}\left( v_{n}\left( \frac{3}{2},-\frac{5}{12}%
\right) -\gamma \right) =\frac{1}{12}.  \tag{2}
\end{equation}%
Now the sequence $\left( v_{n}\left( \frac{3}{2},-\frac{5}{12}\right)
\right) _{n\geq 1}$ has the rate of convergence $n^{-3}$ and the theorem is
proved.$\square $%
\begin{equation*}
\end{equation*}

\textbf{3. Final remarks}\bigskip

In fact, we obtained the sequence%
\begin{equation*}
s_{n}=1+\frac{1}{2}+\cdots +\frac{1}{n-2}+\frac{13}{12\left( n-1\right) }+%
\frac{5}{12n}-\ln n
\end{equation*}%
converging$,$ according to (2), as $n^{-3}$, that is the fastest possible
through all sequences $\left( v_{n}\left( a,b\right) \right) _{n\geq 1}.$

In this case, the best constants $a=\frac{3}{2},$ $b=-\frac{5}{12}$ obtained
in the previous sections can be obtained using another method.

First remark that%
\begin{equation*}
v_{n}\left( a,b\right) =\gamma _{n}+\frac{an+b}{n(n-1)}-\frac{1}{n-1}-\frac{1%
}{n}.
\end{equation*}%
Using the representation of the harmonic sum $h_{n}$ in terms of digamma
function%
\begin{equation*}
h_{n}=\gamma +\frac{1}{n}+\psi \left( n\right) ,
\end{equation*}%
\emph{e.g. }[1, p. 258, Rel. 6.3.2] and the asymptotic formula [1, p. 259,
Rel. 6.3.18]%
\begin{equation*}
\psi \left( z\right) =\ln z-\frac{1}{2z}-\frac{1}{12z^{2}}+\frac{1}{120z^{4}}%
-\frac{1}{252z^{6}}+\cdots ,
\end{equation*}%
we get%
\begin{equation*}
\gamma _{n}=h_{n}-\ln n=\gamma +\frac{1}{n}-\frac{1}{2n}-\frac{1}{12n^{2}}+%
\frac{1}{120n^{4}}-\frac{1}{252n^{6}}+\cdots .
\end{equation*}%
Thus%
\begin{equation*}
v_{n}\left( a,b\right) =\gamma +\frac{(a-\frac{3}{2})n^{2}+\left( b+\frac{5}{%
12}\right) n+\frac{1}{12}}{n^{2}\left( n-1\right) }+\frac{1}{120n^{4}}-\frac{%
1}{252n^{6}}+\cdots .
\end{equation*}%
By analyzing the first fraction in the above representation, the fastest
sequence $\left( v_{n}\left( a,b\right) \right) _{n\geq 1}$ is obtained when
the coefficients $a-\frac{3}{2}$ and $b+\frac{5}{12}$ vanish simultaneously$.
$

On the other hand, let us note that (2) offers us the approximation%
\begin{equation*}
s_{n}-\gamma \approx \frac{1}{12n^{3}},\ \ \ \text{as }n\rightarrow \infty .
\end{equation*}%
We prove the following\bigskip

\textbf{Theorem 2.2. }\emph{For every integer }$n\geq 9,$ \emph{we have}%
\begin{equation*}
\frac{1}{12n^{3}}+\frac{11}{120n^{4}}<s_{n}-\gamma <\frac{1}{12n^{3}}+\frac{%
13}{120n^{4}}.
\end{equation*}

\emph{The left hand side inequality holds for every integer }$n\geq 3.$%
\bigskip

\emph{Proof. }The sequences%
\begin{equation*}
z_{n}=\left( s_{n}-\gamma \right) -\left( \frac{1}{12n^{3}}+\frac{11}{%
120n^{4}}\right) 
\end{equation*}%
and%
\begin{equation*}
t_{n}=\left( s_{n}-\gamma \right) -\left( \frac{1}{12n^{3}}+\frac{13}{%
120n^{4}}\right) 
\end{equation*}%
converges to zero. In order to prove $z_{n}>0$ and $t_{n}<0,$ it suffices to
show that $\left( z_{n}\right) _{n\geq 3}$ is decreasing and $\left(
t_{n}\right) _{n\geq 9}$ is increasing. As%
\begin{equation*}
s_{n+1}-s_{n}=\frac{2}{3n}-\frac{1}{12\left( n-1\right) }+\frac{5}{12\left(
n+1\right) }-\ln \left( 1+\frac{1}{n}\right) ,
\end{equation*}%
we get $z_{n+1}-z_{n}=f\left( n\right) $ and $t_{n+1}-t_{n}=g\left( n\right)
,$ where%
\begin{equation*}
f\left( x\right) =\frac{2}{3x}-\frac{1}{12\left( x-1\right) }+\frac{5}{%
12\left( x+1\right) }-\ln \left( 1+\frac{1}{x}\right) 
\end{equation*}%
\begin{equation*}
-\left( \frac{1}{12\left( x+1\right) ^{3}}+\frac{11}{120\left( x+1\right)
^{4}}\right) +\left( \frac{1}{12x^{3}}+\frac{11}{120x^{4}}\right) 
\end{equation*}%
\begin{equation*}
g\left( x\right) =\frac{2}{3x}-\frac{1}{12\left( x-1\right) }+\frac{5}{%
12\left( x+1\right) }-\ln \left( 1+\frac{1}{x}\right) 
\end{equation*}%
\begin{equation*}
-\left( \frac{1}{12\left( x+1\right) ^{3}}+\frac{13}{120\left( x+1\right)
^{4}}\right) +\left( \frac{1}{12x^{3}}+\frac{13}{120x^{4}}\right) .
\end{equation*}%
Using again Maple software, we obtain%
\begin{equation*}
f^{\prime }\left( x\right) =\frac{P\left( x\right) }{60x^{5}\left(
x-1\right) ^{2}\left( x+1\right) ^{5}}
\end{equation*}%
and%
\begin{equation*}
g^{\prime }\left( x\right) =-\frac{Q\left( x\right) }{60x^{5}\left(
x-1\right) ^{2}\left( x+1\right) ^{5}},
\end{equation*}%
with%
\begin{equation*}
P\left( x\right) =160+1200\left( x-1\right) +2348\left( x-1\right)
^{2}+2055\left( x-1\right) ^{3}+875\left( x-1\right) ^{4}+150\left(
x-1\right) ^{5}
\end{equation*}%
and%
\begin{equation*}
Q\left( x\right) =772\,064+1725\,456\left( x-9\right) +802\,376\left(
x-9\right) ^{2}+164\,805\left( x-9\right) ^{3}
\end{equation*}%
\begin{equation*}
+17\,405\left( x-9\right) ^{4}+930\left( x-9\right) ^{5}+20\left( x-9\right)
^{6}.
\end{equation*}%
Evidently, $f^{\prime }>0$ on $\left( 1,\infty \right) $ and $g^{\prime }<0$
on $\left( 9,\infty \right) .$ It follows that $f$ is strictly increasing on 
$\left( 1,\infty \right) $ and $g$ is strictly decreasing on $\left(
9,\infty \right) .$

As $f\left( \infty \right) =g\left( \infty \right) =0,$ we get $f<0$ on $%
\left( 1,\infty \right) $ and $g>0$ on $\left( 9,\infty \right) .$

It follows that $\left( z_{n}\right) _{n\geq 1}$ is strictly decreasing, $%
\left( t_{n}\right) _{n\geq 9}$ is strictly increasing. As we explained, the
conclusion follows.$\square $\bigskip

\textbf{Acknowledgements. }The work of the second author was supported by a
grant of the Romanian National Authority for Scientific Research,
CNCS-UEFISCDI project number PN-II-ID-PCE-2011-3-0087.%
\begin{equation*}
\end{equation*}

\textbf{References}\bigskip

[1] M. Abramowitz, I. A. Stegun, Handbook of Mathematical Functions with
Formulas, Graphs, and Mathematical Tables, New York: Dover Publications,
1072.

[2] H. Alzer, Inequalities for the gamma and polygamma functions, Abh. Math.
Sem. Univ. Hamburg 68 (1998), 363-372.

[3] G. D. Anderson, R. W. Barnard, M. K. Vamanamurthy, M. Vuorinen,
Inequalities for zero-balanced hypergeometric functions, Trans. Amer. Math.
Soc. 347 (1995), 1713-1723.

[4] C. P. Chen, F. Qi, The best harmonic sequence, arXiv:math/0306233,
available online at: http://arxiv.org/abs/math/0306233.

[5]\ Ch.-P. Chen, Inequalities for the Euler-Mascheroni constant, Appl.
Math. Lett. 23 (2010), 161-164.

[6]\ Ch.-P. Chen, C. Mortici, On a convergence by DeTemple, J. Sci. Arts,
Year 10, No. 2(13) (2010), 271-272.

[7]\ D. W. DeTemple, A quicker convergence to Euler's constant, Amer. Math.
Monthly 100 (5) (1993), 468-470.

[8]\ D. W. DeTemple, A geometric look at sequences that converge to Euler's
constant, College Math. J. 37 (2006), 128-131.

[9]\ E. A. Karatsuba, On the computation of the Euler constant $\gamma $.
Computational methods from rational approximation theory (Wilrijk, 1990),
Numer. Algorithms \ 24 (1-2) (2000), 83-97.

[10]\ K. Knopp, Theory and Applications of Infinite Series, vol 453, Blakie,
London, 1951.

[11]\ C. Mortici, A. Vernescu, An improvement of the convergence speed of
the sequence $\left( \gamma _{n}\right) _{n\geq 1}$ converging to Euler's
constant, An. \c{S}tiin\c{t}. Univ. Ovidius Constan\c{t}a 13 (1) (2005),
97-100.

[12]\ C. Mortici, Product approximation via asymptotic integration, Amer.
Math. Monthly 117 (5) (2010), 434-441.

[13]\ C. Mortici, A. Vernescu, Some new facts in discrete asymptotic
analysis, Math. Balkanika (NS) 21 (Fasc. 3-4) (2007) 301-308.

[14]\ C. Mortici, A Refinement of Chen-Qi Inequality on the Harmonic Sum,
Bul. Univ. Petrol Gaze din Ploie\c{s}ti, Vol LXII No. 1 (2010), 109-112.

[15]\ C. Mortici, Optimizing the rate of convergence in some new classes of
sequences convergent to Euler's constant, Anal. Appl. (Singap.) 8 (1)
(2010), 99-107.

[16]\ C. Mortici, Improved convergence towards generalized Euler-Mascheroni
constant, Appl. Math. Comput. 215 (9) (2010), 3443-3448.

[17]\ C. Mortici, On new sequences converging towards the Euler-Mascheroni
constant, Comput. Math. Appl. 59, 8 (2010), 2610-2614.

[18]\ G. Polya, G. Szeg\H{o}, Problems and Theorems in Analysis, vol I and
II, Springer Verlag, Berlin, Heidelberg, 1972.

[19]\ S. -L. Qiu, M. Vuorinen, Some properties of the gamma and psi
functions with applications, Math. Comp. 74 (250) (2005), 723-742.

[20]\ S. R. Tims, J. A. Tyrrell, Approximate evaluation of Euler's constant,
Math. Gaz. 55 (1971), 65-67.

[21]\ L. Toth, Problem E3432, Amer. Math. Monthly 98 (3) (1991), 264.

[22]\ A. Vernescu, A new accelerate convergence to the constant of Euler,
Gazeta Matem. Ser. A, Bucharest XVII(XCVI) (4) (1999), 273-278.

[23]\ R. M. Young, Euler's constant, Math. Gaz. 75 (1991), 187-190.

\end{document}